\documentclass[11pt,a4paper]{amsart}

\usepackage[T1]{fontenc}
 \usepackage[english]{babel}
 \usepackage[latin1]{inputenc}
\usepackage[dvips]{graphics}
\usepackage{color}
\usepackage{url}
\usepackage{stmaryrd}
\usepackage{rotating}
\usepackage{longtable}
\usepackage{srcltx}

\usepackage{mathabx}
\usepackage{amsmath}
\usepackage{amsthm}
\usepackage{amssymb}
\usepackage{bbm}
\usepackage{enumerate}
\usepackage[all,poly,necula]{xy}

\DeclareSymbolFont{bbdold}{U}{bbold}{m}{n}
\DeclareSymbolFontAlphabet{\mathbbd}{bbdold}

\theoremstyle{plain}
\newtheorem{thm}{Theorem}[section]

\newtheorem{prop}[thm]{Proposition}

\newtheorem{lem}[thm]{Lemma}

\newtheorem{cor}[thm]{Corollary}

\newtheorem*{conjsn}{Conjecture}

\theoremstyle{definition}
\newtheorem{df}[thm]{Definition}

\newtheorem{nt}[thm]{Notation}

\theoremstyle{remark}
\newtheorem*{demo}{Proof}
\newtheorem{rem}[thm]{Remark}

\newtheorem{rems}[thm]{Remarks}

\newcounter{thmnct}

\DeclareMathOperator{\id}{Id}

\DeclareMathOperator{\ima}{im}

\DeclareMathOperator{\red}{_{\!red}}

\DeclareMathOperator{\Gr}{Gr}
\DeclareMathOperator{\rg}{rk}

\DeclareMathOperator{\Hom}{Hom}

\DeclareMathOperator{\md}{mod}

\DeclareMathOperator{\op}{op}

\DeclareMathOperator{\ind}{ind}
\DeclareMathOperator{\irr}{irr}

\DeclareMathOperator{\cyc}{cyc}
\DeclareMathOperator{\rd}{red}
\DeclareMathOperator{\triv}{triv}
\DeclareMathOperator{\Def}{Def}
\DeclareMathOperator{\inW}{in}
\DeclareMathOperator{\injW}{inj}
\DeclareMathOperator{\sym}{sym}
\DeclareMathOperator{\out}{out}
\DeclareMathOperator{\Trop}{Trop}

\newcommand{\lla}{\langle \langle}
\newcommand{\rra}{\rangle \rangle}
\renewcommand{\leq}{\leqslant}
\renewcommand{\geq}{\geqslant}

\newcommand{\Er}{\mathcal{E}}

\renewcommand{\Pr}{\mathcal{P}}
\newcommand{\N}{\mathbb{N}}
\newcommand{\Z}{\mathbb{Z}}
\newcommand{\Q}{\mathbb{Q}}

\newcommand{\C}{\mathbb{C}}

\newcommand{\g}{\mathbf{g}}
\newcommand{\rgr}{\mathbf{r}}
\newcommand{\sgr}{\mathbf{s}}
\newcommand{\egr}{\mathbf{e}}
\newcommand{\h}{\mathbf{h}}

\renewcommand{\epsilon}{\varepsilon}

\newcommand{\inter}{\cap}

\newcommand{\restr}[1]{\!|_{#1}}

\newcommand{\m}{\mathfrak{m}}

\newcommand{\ig}{\mathbf{i}}

\newcommand{\tens}{\otimes}

\newcommand{\cqfd}{\qed}
\renewcommand{\phi}{\varphi}

\renewcommand{\bar}[1]{\overline{#1}}

\newcommand{\comment}[1]{}
\newcommand{\vecv}[2]{\left(\!\!
    \begin{array}{c}
      #1 \\
      #2
    \end{array}
    \!\!\right)}
\newcommand{\vech}[2]{\left(\!\!
    \begin{array}{cc}
      #1 & #2
    \end{array}
    \!\!\right)}
\newcommand{\mat}[4]{\left(\!\!
    \begin{array}{cc}
      #1 & #2 \\
      #3 & #4
    \end{array}
    \!\!\right)}

\newcommand{\inj}{\hookrightarrow}
\newcommand{\surj}{\twoheadrightarrow}

\def\citeb#1#2{\cite[#1]{#2}}

\renewcommand{\tilde}[1]{\widetilde{#1}}

\definecolor{orange}{rgb}{1,0.90,0.5}
\let\oldmarginpar\marginpar
\renewcommand\marginpar[1]{(**)\oldmarginpar[\raggedleft\footnotesize \fcolorbox{blue}{orange}{\parbox{\marginparwidth}{\color{blue}{(**) #1}}}]%
{\raggedright\footnotesize \fcolorbox{blue}{orange}{\parbox{\marginparwidth}{\color{blue}{(**) #1}}}}}

\newcommand{\forloop}[5][1] { \setcounter{#2}{#3} \ifthenelse{#4} { #5 \addtocounter{#2}{#1} \forloop[#1]{#2}{\value{#2}}{#4}{#5} }{ } } 
\newcounter{i}
\newcounter{j}
\newcommand{\dets}[4]{
    \displaystyle \left|\;
    \begin{matrix}
      \forloop{i}{1}{\value{i} < 7}{
       \forloop{j}{1}{\value{j} < 7}{
        \ifthenelse{\value{i} < #1 \or \value{i} > #2 \or \value{j} < #3 \or \value{j} > #4}{\ifthenelse{\value{i}<\value{j}}{\cdot}{\ifthenelse{\value{i} = \value{j}}{1}{}}}{\bullet} \ifthenelse{\value{j}<6}{&}{}
       } \ifthenelse{\value{i}<6}{\\}{}
      }
    \end{matrix}
    \;\right|^{\vphantom{M^{M^M}}}_{\vphantom{M_{M_M}}}}

\begin{document}
\setlength{\parindent}{5mm}
\title[Group species with potentials]{Mutations of group species with potentials and their representations. \\Applications to cluster algebras.}
\author{Laurent Demonet}
\address{Max Planck Institut f\"ur Mathematik, Vivatsgasse 7, D-53111 Bonn}
\email{Laurent.Demonet@normalesup.org}

\date{}

\begin{abstract}
 This article tries to generalize former works of Derksen, Weyman and Zelevinsky about skew-symmetric cluster algebras to the skew-symmetrizable case. We introduce the notion of group species with potentials and their decorated representations. In good cases, we can define mutations of these objects in such a way that these mutations mimic the mutations of seeds defined by Fomin and Zelevinsky for a skew-symmetrizable exchange matrix defined from the group species. These good cases are called non-degenerate. Thus, when an exchange matrix can be associated to a non-degenerate group species with potential, we give an interpretation of the $F$-polynomials and the $\g$-vectors of Fomin and Zelevinsky in terms of the mutation of group species with potentials and their decorated representations. Hence, we can deduce a proof of a serie of combinatorial conjectures of Fomin and Zelevinsky in these cases. Moreover, we give, for certain skew-symmetrizable matrices a proof of the existance of a non-degenerate group species with potential realizing this matrix. On the other hand, we prove that certain skew-symmetrizable matrices can not be realized in this way.
\end{abstract}

\maketitle
\tableofcontents

\section{Introduction}

The aim of this paper is to extend the results of \cite{DeWeZe08} and \cite{DeWeZe} to the case of skew-symmetrizable exchange matrices. Unfortunately, the techniques presented here do not work in any situation, but nevertheless in some important cases.

For this, we introduce \emph{group species with potential} (GSP), which can be seen as quivers with potential with more than one idempotent at each vertex. Thus, we can also define a Jacobian ideal and a Jacobian algebra and study their representations. More precisely, we define the notion of a group species with potential with a \emph{decorated representation} (GSPDR) and the notion of the mutation of a GSPDR at a vertex $k$ (which is called the direction of the mutation). In \emph{good cases}, we can mutate a GSPDR as many times as we want in any direction. In this case, the underlying GSP is called \emph{non-degenerate}. Moreover, we can associate to certain GSPs, called locally free, a skew-symmetrizable matrix in such a way that the mutation we introduce projects, when it exists, to the mutation of matrix introduced by Fomin and Zelevinsky \cite{FoZe02}. Any skew-symmetrizable matrix can be reached in this way using a locally free GSP. The hard problem is to find which skew-symmetrizable matrix can be reached using a non-degenerate GSP. It is the case of matrices of the form $DS$ where $D$ is diagonal with positive integer coefficients and $S$ is skew-symmetric with integer coefficients. It is also the case for the skew-symmetrizable matrices which occur in the situation of \cite{De}, in particular in all acyclic cases. Nevertheless, it is not always true, as shown by the counterexample at the end of section \ref{exemple}. The techniques presented in \cite{DeWeZe08} work here almost in the same way. The only problem is that it is not always the case that for any $2$-cycle, there exists a potential canceling it (this fact is very easy in the context of \cite{DeWeZe08}).

We now explain the content of this article in more details. Let $K$ be an algebraically closed field. Let $I$ be a finite set and $E = \bigoplus_{i \in I} K[\Gamma_i]$ where, for each $i$, $\Gamma_i$ is a finite group whose cardinal is not divisible by the characteristic of $K$. Let also $A$ be an $(E, E)$-bimodule. This data is called a \emph{group species} and its \emph{complete path algebra} is
 $$E \lla A \rra = \prod_{n \in \N} A^{\tens n}.$$
A potential $S$ on this group species can be seen as a (maybe infinite) linear combination of cyclic path, up to rotation. It permits to construct a two sided ideal $J(S)$, called the \emph{Jacobian ideal} and a quotient algebra $\Pr(A, S) = E \lla A \rra / J(S)$ called the \emph{Jacobian algebra}. A \emph{decorated representation} of the GSP is a pair consisting of a $\Pr(A, S)$-module $X$ and an $E$-module $V$. In sections \ref{mut} and \ref{mutrep}, we define the mutation of a GSP with a decorated representation (GSPDR). This mutation is well defined if the group species has no loop and is $2$-acyclic (that is, for any $i \in I$, $E_i (A \oplus A \tens_E A) E_i = 0$, where $E_i = K[\Gamma_i] \subset E$). 
 
In what follows, we suppose that the $\Gamma_i$ are commutative and that the GSP is \emph{locally free}, that is, for any $i, j \in I$, $E_i A E_j$ is a free $E_i$-left module and a free $E_j$-right module. In section \ref{exc}, we define the exchange matrix $B$ of a the group species by
 $$b_{ij} = \dim_{E_j} A_{ji} - \dim_{E_j} A_{ij}^*.$$
 Thus, the mutation of GSPDRs descends to the mutation of matrices defined by Fomin and Zelevinsky \cite{FoZe02}. In section \ref{nondeg}, we discuss a class of matrices, namely those of the form $DS$, for which there is always a non-degenerate GSP. Moreover, we remark that there exists also non-degenerate GSP in all cases which are categorified in \cite{De} (because the endomorphisms rings of cluster-tilting objects constructed in \cite{De} are Jacobian algebras). Remark also that there is no chance, with definitions given here, to construct non-degenerate GSPs for any skew-symmetrizable matrix, as shown by the counterexample ending section \ref{exemple}.

Following the ideas of \cite{DeWeZe}, we explain in section \ref{fpol} how to reinterpret the $F$-polynomials and $\g$-vectors defined in \cite{FoZe07} in terms of GSPDRs and their mutations. We deduce in section \ref{clust} that, when a skew-symmetrizable matrix can be obtained from a non-degenerate GSP, then the following conjectures are true:

\begin{conjsn}[\citeb{conjecture 5.4}{FoZe07}]
 For any $\ig \in I^n$ and $k \in I$, $F^B_{k;\ig}$ has constant term $1$.
\end{conjsn}

\begin{conjsn}[\citeb{conjecture 5.5}{FoZe07}]
 For any $\ig \in I^n$ and $k \in I$, $F^B_{k;\ig}$ has a maximum monomial for divisibility order with coefficient $1$.
\end{conjsn}

\begin{conjsn}[\citeb{conjecture 7.12}{FoZe07}]
 For any $\ig \in I^n$, $k \in I$, we denote by $k\ig$ the concatenation of $(k)$ and $\ig$. Let $j \in I$ and $(g_i)_{i \in I} = \g^B_{j; \ig}$ and $(g'_i)_{i \in I} = \g^{\mu_k(B)}_{j; k\ig}$. Then we have, for any $i \in I$,
 $$g'_i = \left\{\begin{array}{ll}
  -g_i & \quad \text{if } i = k; \\
  g_i + \max(0,b_{ik}) g_k - b_{jk} \min(g_k, 0) & \quad \text{if } i \neq k.
  \end{array} \right.$$
\end{conjsn}

\begin{conjsn}[\citeb{conjecture 6.13}{FoZe07}]
 For any $\ig \in I^n$, the vectors $\g^B_{i;\ig}$ for $i \in I$ are sign-coherent. In other terms, for $i,i', j \in I$, the $j$-th components of $\g^B_{i;\ig}$ and $\g^B_{i';\ig}$ have the same sign.
\end{conjsn}

\begin{conjsn}[\citeb{conjecture 7.10(2)}{FoZe07}]
 For any $\ig \in I^n$, the vectors $\g^B_{i;\ig}$ for $i \in I$ form a $\Z$-basis of $\Z^I$.
\end{conjsn}

\begin{conjsn}[\citeb{conjecture 7.10(1)}{FoZe07}]
 For any $\ig, \ig' \in I^n$, if we have
 $$\sum_{i \in I} a_i \g_{i; \ig}^B = \sum_{i \in I} a'_i \g_{i; \ig'}^B$$
 for some nonnegative integers $(a_i)_{i \in I}$ and $(a'_i)_{i \in I}$, then there is a permutation $\sigma \in \mathfrak{S}_I$ such that
 for every $i \in I$, 
 $$a_i = a'_{\sigma(i)} \quad \text{and} \quad a_i \neq 0 \Rightarrow \g_{i;\ig}^B = \g_{\sigma(i);\ig'}^B \quad \text{and} \quad a_i \neq 0 \Rightarrow F_{i;\ig}^B = F_{\sigma(i);\ig'}^B.$$
 In particular, $F_{i; \ig}^B$ is determined by $\g_{i; \ig}^B$.
\end{conjsn}

Thus, as stated in \cite[remark 7.11]{FoZe07}, if $B$ is a full rank skew-sym\-me\-tri\-zable matrix which correspond to a non-degenerate GSP, then the cluster monomials of a cluster algebra with exchange matrix $B$ are linearly independent.

\section{Group species and path algebras}

Let $K$ be a field.

\begin{df}
 A {\em group species} is a triple $(I, (\Gamma_i)_{i \in I}, (A_{ij})_{(i,j) \in I^2})$ where $I$ is a finite set, for each $i \in I$, $\Gamma_i$ is a finite group and for each $(i,j) \in I^2$, $A_{ij}$ is a finite dimensional $(K[\Gamma_i], K[\Gamma_j])$-bimodule (the first acting on the left and the second on the right).
\end{df}

Fix now such a group species $Q = (I, (\Gamma_i)_{i \in I}, (A_{ij})_{(i,j) \in I^2})$ 

\begin{df}
 A {\em representation of $Q$} is a pair $((V_i)_{i \in I}, (x_{ij})_{(i,j) \in I^2})$ where for each $i \in I$, $V_i$ is a right finite dimensional $K[\Gamma_i]$-module and for each $(i,j) \in I^2$, 
 $$x_{ij} \in \Hom_{\Gamma_j} (V_i \tens_{\Gamma_i} A_{ij}, V_j).$$
\end{df}

\begin{df}
 Let $((V_i)_{i \in I}, (x_{ij})_{(i,j) \in I^2})$ and $((V'_i)_{i \in I}, (x'_{ij})_{(i,j) \in I^2})$ be two representations of $Q$. A {\em morphism} from the first one to the second one is a family $(f_i)_{i \in I} \in \prod_{i \in I} \Hom_{\Gamma_i}(V_i, V'_i)$ such that for each $(i, j) \in I^2$ the following diagram commute :
 $$\xymatrix{
  V_i \tens_{\Gamma_i} A_{ij} \ar[d]_{f_i \tens \id_{A_{ij}}} \ar[r]^-{x_{ij}} & V_j \ar[d]^{f_j} \\
  V'_i \tens_{\Gamma_i} A_{ij} \ar[r]_-{x'_{ij}} & V'_j 
 }$$
\end{df}

\begin{rems}
 \begin{itemize}
  \item The previous definitions give rise to an abelian category.
  \item If for each $i \in I$, $\Gamma_i$ is the trivial group, we get back the classical definition of a quiver (up to the choice of a basis of each $A_{ij}$) and of the category of representations of a quiver.
  \item If for each $i \in I$, $K[\Gamma_i]$ is replaced by a division algebra, we obtain the usual definition of a species (see for example \cite{DlRi76}).
 \end{itemize} 
\end{rems}

\begin{df}
 For each $i \in I$, denote $E_i = K[\Gamma_i]$. Denote also $E = \bigoplus_{i \in I} E_i$ and $A = \bigoplus_{(i, j) \in I^2} A_{ij}$. Thus, we put the natural $(E, E)$-bimodule structure on $A$ and define the graded algebras
 $$E \langle A \rangle = \bigoplus_{n \in \N} A^{\tens n} \quad \text{and} \quad E \lla A \rra = \prod_{n \in \N} A^{\tens n}$$
 the first one being called the \emph{path algebra} of the group species and the second one the \emph{complete path algebra} of the group species (note that every tensor product is taken over $E$).
\end{df}

\begin{rems} \label{dfm}
 \begin{itemize}
  \item As usual for quiver, the category of representations of a group species is equivalent to the category of finite dimensional right modules over its path algebra. Moreover, the category of nilpotent representations of a group species is equivalent to the category of finite dimensional right modules over its complete path algebra.
  \item If one denotes 
    $$\m = \prod_{n > 0} A^{\tens n} \subset E \lla A \rra$$
    which is clearly a two-sided ideal, then $E \lla A \rra$ becomes a topological algebra for the $\m$-adic topology and $E \langle A \rangle$ is a dense subalgebra of it.
 \end{itemize}
\end{rems}

As in \cite{DeWeZe08}, $\m$ is the unique maximal two-sided ideal of $E \lla A \rra$ not intersecting $E$. Moreover, if we have another group species with the same vertices whose arrows are encoded in the $(E, E)$-bimodule $A'$, then, again as in \cite{DeWeZe08}, the morphisms $\phi$ from $E\lla A \rra$ to $E \lla A' \rra$ such that $\phi \restr{E} =\id_E$ (later called $E$-morphisms) are parameterized in an obvious way by a pair $(\phi^{(1)}, \phi^{(2)})$ where $\phi^{(1)} : A \rightarrow A'$ and $\phi^{(2)} : A \rightarrow \m'^2$ are $(E,E)$-bimodule morphisms. Thus, $\phi$ is an isomorphism if and only if $\phi^{(1)}$ is an isomorphism. Introduce now the analogous of \cite[definition 2.5]{DeWeZe08}:

\begin{df}
 An $E$-automorphism $\phi$ of $E \lla A \rra$ will be called a \emph{change of arrows} if $\phi^{(2)} = 0$ and a \emph{unitriangular automorphism} if $\phi^{(1)} = \id_A$.
\end{df}

Finally, introduce the following useful notation:

\begin{nt}
 For all $i, j \in I$, 
 $$E \langle A \rangle_{ij} = E_i E \langle A \rangle E_j \quad \text{and} \quad E \lla A \rra_{ij} = E_i E \lla A \rra E_j$$
 and for $n \in \N$,
 $$A_{ij}^{\tens n} = A^{\tens n} \inter E \langle A \rangle_{ij} = A^{\tens n} \inter E \lla A \rra_{ij}$$
 so that
 $$E \langle A \rangle_{ij} = \bigoplus_{n \in \N} A_{ij}^{\tens n} \quad \text{and} \quad E \lla A \rra_{ij} = \prod_{n \in \N} A_{ij}^{\tens n}.$$
\end{nt}

\section{Potential and their Jacobian ideals}

Following \cite{DeWeZe08} define:

\begin{df}
 Define
  $$E \lla A \rra_{\cyc} = \frac{E \lla A \rra}{\left[E \lla A \rra, E \lla A \rra\right]}$$
 whose elements are called \emph{potentials} (here, $\left[E \lla A \rra, E \lla A \rra\right]$ is the closure of the two-sided ideal generated by commutators). As $\left[E \lla A \rra, E \lla A \rra\right]$ is generated by its homogeneous elements, we can decompose $E \lla A \rra_{\cyc} = \prod_{n \in \N} A^{\tens n}_{\cyc}$ where
 $$A^{\tens n}_{\cyc} = \frac{A^{\tens n}}{\left[E \lla A \rra, E \lla A \rra\right] \inter A^{\tens n}}$$
 and, if $S \in E \lla A \rra_{\cyc}$, we write $S^{(n)}$ its summand which lies in $A^{\tens n}_{\cyc}$.
\end{df}

\begin{df}
 Define the continuous linear map
 $$\partial : \left(E \lla A \rra\right)^* \tens_k E \lla A \rra \rightarrow E \lla A \rra$$
 in the following way. First remark that $\left(E \lla A \rra\right)^* \simeq \bigoplus_{n \in \N} \left( A^{\tens n}\right)^*$. Then, if $\xi \in \left( A^{\tens n}\right)^*$ and $a_1, a_2, \dots, a_\ell \in A$ define $\partial_\xi(a_1 a_2 \dots a_\ell) = 0$ if $\ell < n$ and 
 $$\partial_\xi(a_1 a_2 \dots a_\ell) = \sum_{j = 1}^\ell \sum_{g, h \in \mathcal{B}} \xi\left(g^{-1} a_j a_{j+1} \dots a_{j+n-1} h\right) h^{-1} a_{j+n} a_{j+n+1} \dots a_{j-1} g$$
 if $\ell \geq n$ where all indices are taken modulo $\ell$ and $\mathcal{B} = \bigcup_{i \in I} \Gamma_i \subset E$. It is easy to see that $\partial$ is well defined and moreover that it vanishes on commutators. Thus, we can descend $\partial$ to a continuous linear map
 $$\partial : \left(E \lla A \rra\right)^* \tens_k E \lla A \rra_{\cyc} \rightarrow E \lla A \rra.$$
\end{df}

\begin{rem}
 With the natural structure of $(E, E)$-bimodule on $\left(E \lla A \rra\right)^*$, one gets, for any $S \in E \lla A \rra_{\cyc}$, that $\xi \mapsto \partial_\xi S$ is a morphism of $(E, E)$-bimodules.
\end{rem} 

\begin{df}
 For a potential $S \in E \lla A \rra_{\cyc}$, define the \emph{Jacobian ideal} $J(S)$ to be the closure of the two-sided ideal of $E \lla A \rra$ generated by the $\partial_\xi(S)$ for $\xi \in A^*$. The quotient $E \lla A \rra / J(S)$ is called the \emph{Jacobian algebra} and is denoted by $\Pr(A, S)$ (we do not keep trace of $(I, (\Gamma_i))$ in the notation because it will be fixed).
\end{df}

Note that every $E$-morphism $\phi: E \lla A \rra \rightarrow E \lla A' \rra$ descends to $\phi: E \lla A \rra_{\cyc} \rightarrow E \lla A' \rra_{\cyc}$.

It is now easy to adapt the proof of \citeb{proposition 3.7}{DeWeZe08}:

\begin{prop}
 Let $S \in E\lla A \rra_{\cyc}$. Every $E$-isomorphism $\phi: E\lla A \rra \rightarrow E\lla A' \rra$ maps $J(S)$ to $J(\phi(S))$ and therefore induces an isomorphism $$\Pr(A, S) \rightarrow \Pr(A', \phi(S)).$$
\end{prop}

\section{Group species with potentials}

For the rest of this article, the data $(I, (\Gamma_i))$ and so $E$ will be fixed. Following the ideas of \cite{DeWeZe08}, define:

\begin{df}
 As before, $A$ is an $(E,E)$-bimodule and we take $S \in E \lla A \rra_{\cyc}$. We say that $(A, S)$ is a group species with potential (GSP for short) if the species has no loop (for all $i \in I$, $E_i A E_i = \{0\}$) and $S \in \prod_{n > 1} A^{\tens n}_{\cyc}$.
\end{df}

\begin{df} \label{reqg}
 Let $(A, S)$ and $(A', S')$ be two GSPs. One says that an $E$-isomorphism $\phi: E\lla A \rra \rightarrow E\lla A' \rra$ is a right-equivalence if $\phi(S) = S'$.
\end{df}

Note that this definition induces a equivalence relation. Moreover, a right equivalence $(A, S) \simeq (A', S')$ induces isomorphisms of $(E, E)$-bimodules $A \simeq A'$, $J(S) \simeq J(S')$ and $\Pr(A, S) \simeq \Pr(A', S')$ as said before.

\begin{nt}
 If $(A, S)$ and $(A', S')$ are two GSPs, define $(A, S) \oplus (A', S') = (A \oplus A', S + S')$ so that $\Pr((A, S) \oplus (A', S'))$ is the completion of $\Pr(A, S) \oplus \Pr(A', S')$ for the product topology. 
\end{nt}

\begin{df}
 We say that a GSP $(A, S)$ is \emph{trivial} if $S \in A^{\tens 2}_{\cyc}$ and $\{\partial_\xi(S)\,|\, \xi \in A^*\} = A$, or, equivalently, if $\Pr(A, S) = E$.
\end{df}

The following easy proposition is an adaptation of \cite[proposition 4.4]{DeWeZe08}:

\begin{prop}
 A GSP $(A, S)$ is trivial if and only if there exist an $(E,E)$-bimodule $B$ and an $(E,E)$-bimodules isomorphism $\phi : A \rightarrow B \oplus B^*$ such that 
 $$\phi(S) = \sum_{b \in \mathcal{B}} b \tens b^*$$
 where $\phi$ is naturally extended to an isomorphism $E \lla A \rra_{\cyc} \rightarrow E \lla B \oplus B^* \rra_{\cyc}$ and the right member does not depend of the choice of a basis $\mathcal{B}$ of $B$.
\end{prop}

One gets also this proposition, similar to \citeb{proposition 4.5}{DeWeZe08}:

\begin{prop} \label{eqenltriv}
 If $(A, S)$ is a GSP and $(B, T)$ is a trivial GSP, then the canonical embedding $E \lla A \rra \inj E \lla A \oplus B \rra$ induces an isomorphism $\Pr(A, S) \simeq \Pr(A \oplus B, S + T)$.
\end{prop}

For a GSP $(A, S)$, we define the \emph{trivial} and \emph{reduced} part of $A$ as the $(E,E)$-bimodules
 $$A_{\triv} = \{\partial_\xi S^{(2)}\,|\, \xi \in A^*\} \quad \text{and} \quad A_{\rd} = A/A_{\triv}.$$
Moreover, we say that $(A, S)$ is reduced if $S^{(2)} = 0$, or, equivalently, if $A_{\triv} = \{0\}$.

Again, the proof of \cite[theorem 4.6]{DeWeZe08} is easy to adapt:

\begin{thm} \label{reqredtriv}
 For any GSP $(A, S)$, there exist $S_{\triv} \in E \lla A_{\triv} \rra$ and  $S_{\rd} \in E \lla A_{\rd} \rra$ such that $(A, S)$ is right equivalent to $(A_{\triv}, S_{\triv}) \oplus (A_{\rd}, S_{\rd})$.

 Moreover, the right equivalence classes of $(A_{\triv}, S_{\triv})$ and $(A_{\rd}, S_{\rd})$ are uniquely determined by the right equivalence class of $(A, S)$.
\end{thm}

\begin{df}
 A group species $(I, (\Gamma_i), A)$ is called \emph{$2$-acyclic} if, for any $i \in I$, $E_i A^{\tens 2} E_i = \{0\}$.
\end{df}

We will see now how to find, as in \cite{DeWeZe08}, algebraic conditions guaranteeing the $2$-acyclicity of the reduced part of a group species. Let $K\left[E \lla A \rra_{\cyc}\right]$ be the ring of polynomial functions on $E \lla A \rra_{\cyc}$ vanishing on all but a finite number of the $A^{\tens n}_{\cyc}$.

For each $S \in E \lla A \rra_{\cyc}$ and $i, j \in I$, define the bilinear form $\alpha_{S,ij}$ by:
 \begin{align*} 
  A_{ij}^* \times A_{ji}^* &\rightarrow K \\
       (f, g) & \mapsto \mathop{\sum_{\gamma \in \Gamma_i}}_{\gamma' \in \Gamma_j} \left[(\gamma' f \gamma^{-1} \tens \gamma g \gamma'^{-1})\left(S^{(2)}\right) + (\gamma g \gamma'^{-1} \gamma' f \gamma^{-1})\left(S^{(2)}\right)\right].
 \end{align*}

First, an easy lemma:

\begin{lem}
 Let $i, j \in I$. The followings are equivalent:
 \begin{enumerate}
  \item there exists $S \in E \lla A \rra_{\cyc}$ such that $\alpha_{S,ij}$ is of maximal rank;
  \item either $A^*_{ij}$ is a subbimodule of $A_{ji}$ or $A^*_{ji}$ is a subbimodule of $A_{ij}$.
 \end{enumerate}
\end{lem}

\begin{demo}
 We clearly have $\alpha_{S,ij} = \alpha_{S,ji}$ for any $S$ and therefore, one can suppose without loss of generality that $\dim_K A_{ij} \leq \dim_K A_{ji}$. Suppose that $\alpha_{S, ij}$ is of maximal rank. In any basis, the matrix of $\alpha_{S, ij}$ is the matrix of  $A^*_{ij} \rightarrow A_{ji} : \xi \mapsto \partial_{\xi}(S^{(2)})$ and therefore, $A^*_{ij}$ is a subbimodule of $A_{ji}$.

 Reciprocally, suppose that $A^*_{ij}$ is a subbimodule of $A_{ji}$. Thus, if $\mathcal{B}$ is a basis of $A_{ij}$, define
 $$S = \sum_{a \in \mathcal{B}} a \tens a^*$$
 where $a^* \in A^*_{ij}$ is identified with its image in $A_{ji}$. Then, it is clear that $\alpha_{S, ij}$ is of maximal rank. \cqfd
\end{demo}

Again, it is easy to generalize \citeb{proposition 4.15}{DeWeZe08}:

\begin{prop} \label{eff2c}
 The reduced part of a GSP $(A, S)$ is $2$-acyclic if and only if, for any $i, j \in I$, $\alpha_{S, ij}$ is of maximal rank. This condition is open. Moreover, if, for any $i, j \in I$, either $A_{ij}^*$ is a subbimodule of $A_{ji}$, either $A_{ji}^*$ is a subbimodule of $A_{ij}$, then there is a non empty Zariski open subset $U$ of $E \lla A \rra_{\cyc}$, a $2$-acyclic $(E,E)$-bimodule $A'$ and a regular map $H: U \rightarrow E \lla A' \rra_{\cyc}$ such that for any $S \in U$, $(A_{\rd}, S_{\rd})$ is right equivalent to $(A', H(S))$.
\end{prop}         

\begin{demo}
 The arguments are the same than in \cite{DeWeZe08}. For each $i, j \in I^2$, choose $\bar A_{ij}^* \subset A_{ij}^*$ such that $\bar A_{ij}^* = A_{ij}^*$ or $\bar A_{ij}^* \simeq A_{ji}$. Let $U$ to be the non-empty open subset of $E \lla A \rra_{\cyc}$ containing the $S$ such that for all $i,j \in I$, $\alpha_{S,ij}\restr{\bar A_{ij}^* \times \bar A_{ji}^*}$ is non-degenerate (it corresponds to the non-vanishing of a fixed maximal minor of $\alpha_{S, ij}$). Define $A'$ to be the intersection of the kernels of the elements of the $\bar A_{ij}^*$. Then the construction of $H$ follows the proof of \cite[theorem 4.6]{DeWeZe08}. \cqfd
\end{demo}

\section{Mutations of group species with potential} \label{mut}

Let $(A, S)$ and $k \in I$ be a vertex such that $E_k A^{\tens 2} E_k = \{0\}$ (we say that $(A, S)$ is $2$-acyclic at $k$). We suppose also that for any $i \in I$, the characteristic of $K$ does not divide $\# \Gamma_i$. As in \cite[\S 5]{DeWeZe08}, one defines $\tilde \mu_k(A, S) = (\tilde A, \tilde S)$ where, if $i, j \in I$, 
$$\tilde A_{ij} = \left\{ \begin{array}{ll}
 A_{ji}^* & \quad \text{if } k \in \{i, j\}; \\
 A_{ij} \oplus A_{ik} \tens_{E_k} A_{kj} & \quad \text{otherwise}.
\end{array} \right.$$
In other terms,
$$\tilde A = \bar E_k A \bar E_k \oplus A E_k A \oplus (E_k A)^* \oplus (A E_k)^*$$
where $\bar E_k = \bigoplus_{i \neq k} E_i$. Let now $[-] : \bar E_k E \lla A \rra \bar E_k \rightarrow E \lla \tilde A \rra$ be the morphism of $k$-algebras generated by $[a] = a$ if $a \in \bar E_k A \bar E_k$ and $[ab] = ab \in A E_k A$ if $a \in A E_k$ and $b \in E_k A$ which is well defined because $(A, S)$ has no loop. Again, because $(A, S)$ has no loop, every potential $S \in E \lla A \rra_{\cyc}$ has a representative in $\bar E_k E \lla A \rra \bar E_k$ and it is easy to see that $[-]$ descends to a map 
 $$[-]: E \lla A \rra_{\cyc} \rightarrow E \lla \tilde A \rra_{\cyc}.$$
Moreover, as for any $i \in I$ the characteristic of $K$ does not divide $\# \Gamma_i$, we have a canonical sequence of isomorphisms
 \begin{align*} \Hom_E\left( A E_k A, A E_k A \right) &\simeq (A E_k A)^* \tens_E A E_k A \simeq (A E_k \tens_E E_k A)^* \tens_E A E_k A \\ &\simeq  (E_k A)^* \tens_E (A E_k)^* \tens_E A E_k A \subset E \lla \tilde A \rra \end{align*}
and we define $\Delta_k(A)$ to be the image of $\id_{A E_k A}$ through this isomorphism. Thus, define
 $$\tilde S = [S] + \Delta_k(A).$$

The proof of \citeb{proposition 5.1}{DeWeZe08} can be easily generalized:

\begin{prop}
 If $(A', S')$ is another GSP such that $E_k A' = A' E_k = \{0\}$, then 
  $$\tilde \mu_k(A \oplus A', S + S') = \mu_k(A, S) \oplus (A', S').$$
\end{prop}

Now, the proof of \citeb{theorem 5.2}{DeWeZe08} is easy to generalize:

\begin{thm}
 \label{recm} The right-equivalence class of the GSP $\tilde \mu_k(A, S)$ is fully determined by the right-equivalence class of $(A, S)$.
\end{thm}

\begin{df}
 Using theorem \ref{recm} together with theorem \ref{reqredtriv}, the right-equivalence class of the reduced part of $\tilde \mu_k(A,S)$ is fully determined by the right-equivalence class of $(A, S)$. Thus we can define the map $\mu_k$ from the set of right-equivalence classes which are $2$-acyclic at $k$ to itself. It is called the \emph{mutation at vertex $k$}.
\end{df}

Again, the proof of \cite[theorem 5.7]{DeWeZe08} is easy to generalize:

\begin{thm}
 $\mu_k$ is an involution.
\end{thm}

Let us also remark that \citeb{proposition 6.1}{DeWeZe08}, \citeb{proposition 6.4}{DeWeZe08} and \citeb{corollary 6.6}{DeWeZe08} can be generalized:

\begin{prop} \label{enlevevertex}
 The algebras $\bar E_k \Pr(A, S) \bar E_k$ and $\bar E_k \Pr\left(\tilde \mu_k(A, S)\right) \bar E_k$ are isomorphic.
\end{prop}

\begin{prop}
 The Jacobian algebra $\Pr(A, S)$ is finite-dimensional if and only if $\Pr\left(\tilde \mu_k(A, S)\right)$ is.
\end{prop}

\begin{cor}
 The Jacobian algebras $\bar E_k \Pr(A, S) \bar E_k$ and $\bar E_k \Pr\left(\mu_k(A, S)\right) \bar E_k$ are isomorphic and $\Pr(A, S)$ is finite-dimensional if and only if $\Pr\left(\mu_k(A, S)\right)$ is.
\end{cor}

As stated in \cite[remark 6.8]{DeWeZe08}, the following definition makes sense:

\begin{df}
 We define the \emph{deformation space of} $(A, S)$ to be
 $$\Def(A, S) = \frac{\Pr(A, S)}{E + \left[\Pr(A, S), \Pr(A, S)\right]}$$
 where $\left[\Pr(A, S), \Pr(A, S)\right]$ is the closure of the two-sided ideal of $\Pr(A,S)$ generated by the commutators.
\end{df}

Thus, let us introduce the following extension of \citeb{proposition 6.9}{DeWeZe08}:

\begin{prop}
 We have an isomorphism:
 $$\Def(A,S) \simeq \Def\left(\tilde \mu_k(A, S)\right).$$
\end{prop}

\begin{demo}
 It is enough to prove that
 $$\frac{\bar E_k \Pr(A, S) \bar E_k}{\bar E_k + \left[\bar E_k \Pr(A, S) \bar E_k, \bar E_k \Pr(A, S) \bar E_k\right]} \inj \Def(A, S)$$
 is in fact an isomorphism (which is true because $A$ has no loop) and to use proposition \ref{enlevevertex}.
\end{demo}

As in \cite{DeWeZe08},

\begin{df}
 The GSP $(A, S)$ is called \emph{rigid} if $\Def(A, S) = \{0\}$.
\end{df}

\begin{cor}
 The GSP $(A, S)$ is rigid if and only if $\mu_k(A, S)$ is.
\end{cor}

\section{Exchange matrices} \label{exc} 

We suppose now that $A$ has neither loop nor $2$-cycle (that is $A^{\tens 1}_{\cyc} = A^{\tens 2}_{\cyc} = \{0\}$). We suppose also that for any $(i, j) \in I^2$, $A_{ij}$ is a free left $E_i$-module and a free right $E_j$-module (we will call it a \emph{locally free} GSP). Define the matrix $B = B(A) = B(A, S)$ to be the matrix with rows and columns indexed by $I$ and coefficients
 $$b_{ij} = \dim_{E_j} A_{ji} - \dim_{E_j} A_{ij}^*$$
 (by default, dimension are taken relatively to the left module structure).
This matrix is clearly skew-symmetrizable since
 $$\# \Gamma_j \times b_{ij} = \dim_K A_{ji} - \dim_K A_{ij}^*.$$

\begin{df}
 The matrix $B$ is called the \emph{exchange matrix} of $A$.
\end{df}

The following proposition justifies this generalization of \cite{DeWeZe08}:
\begin{prop} \label{consgsp}
 Every skew-symmetrizable matrix $B$ can be reached in this way from a GSP.
\end{prop}

\begin{demo}
 Let $B$ be a skew-symmetrizable matrix and $D = (d_i)_{i \in I}$ be a diagonal matrix with positive integer coefficients such that $BD$ is skew-symmetric. Let $\Gamma_i = \Z / d_i \Z$ and for $(i, j) \in I^2$ such that $b_{ij} > 0$,
 $$A_{ji} = K\left[\Z / (d_j b_{ij}) \Z\right] = K\left[\Z / (-d_i b_{ji}) \Z\right]$$
 which is a left and right free $(\Gamma_j, \Gamma_i)$-bimodule. It is clear that $A = \bigoplus_{i,j \in I} A_{ij}$ has exchange matrix $B$. \cqfd
\end{demo}

\begin{prop}
 Let $k \in I$.
 \begin{enumerate}
  \item The GSP $\tilde \mu_k(A, S)$ is locally free.
  \item If $\mu_k(A, S)$ is $2$-acyclic then it is locally free.
  \item If $\mu_k(A, S)$ is $2$-acyclic then
   $$\mu_k(B(A, S)) = B(\mu_k(A, S))$$
   where the $\mu_k$ on the left hand is the one defined in \cite{FoZe02}. Namely:
   $$b'_{ij} = \left\{ \begin{array}{ll}
    \displaystyle -b_{ij} & \quad \text{if } k \in \{i, j\} \\
    \displaystyle b_{ij} + \frac{b_{ik}\,|b_{kj}| + |b_{ik}|\,b_{kj}}{2} & \quad \text{otherwise}
    \end{array} \right.$$
   if $B' = \mu_k(B)$.
 \end{enumerate}
\end{prop}

\begin{demo}
 \begin{enumerate}
  \item First of all, it is clear that for $i \in I$, $E_i^* \simeq E_i$ as $(E_i, E_i)$-bimodules (as $E_i$ is finite dimensional). Thus, for any $i$, $A_{ik}^*$ and $A_{ki}^*$ are left and right free modules. Moreover, as a right module, 
   $$A_{ik} \tens_{E_k} A_{kj} \simeq A_{kj}^{\dim_{E_k} \left(A_{ik}^*\right)}$$
   and, as a left module,
   $$A_{ik} \tens_{E_k} A_{kj} \simeq A_{ik}^{\dim_{E_k} \left(A_{kj}\right)}$$
   which ends the proof that $\tilde \mu_k(A, S)$ is locally free.
  \item If one denotes $(\tilde A, \tilde S) = \tilde \mu_k(A, S)$, one has
   $$\tilde A = \tilde A_{\rd} \oplus \tilde A_{\triv}$$
   As $\tilde A_{\rd}$ is $2$-acyclic, for any $i, j \in I$, $\tilde A_{\rd, ij} = 0$ or $\tilde A_{\rd, ji} = 0$. Suppose that $\tilde A_{\rd, ij} = 0$. Hence $\tilde A_{\triv, ji} \simeq \tilde A_{\triv, ij}^* \simeq \tilde A_{ij}^*$ is left and right free (thanks to the previous point). Moreover, $\tilde A_{ji} = \tilde A_{\rd, ji} \oplus \tilde A_{\triv, ji}$ and, as the categories of left $E_j$-modules and right $E_i$-modules are Krull-Schmidt, $\tilde A_{\rd, ji}$ is left and right free.
  \item It is enough to remark that
   $$\dim_{E_i} A_{ik} \tens_{E_k} A_{kj} = \dim_{E_i} A_{ik}^{\dim_{E_k} A_{kj}} = \dim_{E_i}(A_{ik})\dim_{E_k}(A_{kj})$$ 
  and that
   $$\dim_{E_i} (A_{jk} \tens_{E_k} A_{ki})^* = \dim_{E_i} \left(A_{ki}^*\right)^{\dim_{E_k} A_{jk}^*} = \dim_{E_i}(A_{ki}^*) \dim_{E_k} (A_{jk}^*)$$ 
  and to use the definition and the duality $A_{\triv, ij} \simeq A_{\triv, ji}^*$. \cqfd
 \end{enumerate}
\end{demo}

\begin{df}
 The group species is said to be \emph{globally free} if, for any $i, j \in I$, $A_{ij}$ is a free $(E_i, E_j)$-bimodule (that is a free $E_i \tens_K E_j^{\op}$-module). 
\end{df}

\begin{rem}
 The class of globally free group species is stable under mutation.
\end{rem}

\begin{prop}
 If a matrix is of the form $DB$, where $D$ is diagonal with positive integer coefficients and $B$ is skew-symmetric, then the group species constructed in proposition \ref{consgsp} is globally free.
\end{prop}

\section{Existance of nondegenerate potentials} \label{nondeg}

If $(I, (\Gamma_i), A)$ is a group species without loop nor $2$-cycle, a potential $S \in E \lla A \rra_{\cyc}$ will be said to be \emph{non-degenerate} if every sequence of mutation going from $(A, S)$ yields to a $2$-acyclic GSP. 

We cite the following adapted result, whose proof is the same than the proof of \citeb{corollary 7.4}{DeWeZe08}:

\begin{thm}
 If the group species is globally free then there is a countable number of non-constant polynomials in $K\left[E \lla A \rra_{\cyc}\right]$ such that the non-vanishing of these polynomials on $S \in E \lla A \rra_{\cyc}$ implies that $S$ is non-degenerate. In particular if $K$ is uncountable, there exist non-degenerate potentials.
\end{thm}

\begin{demo}
 The only thing to change is that, if the group species is globally free, then for each $i, j \in I$, either $A_{ij}^*$ is a subbimodule of $A_{ji}$, or $A_{ji}^*$ is a subbimodule of $A_{ij}$ and, therefore, proposition \ref{eff2c} can be applied. \cqfd
\end{demo}

\begin{rem}
 It is also easy to prove that for any skew-symmetrizable matrix $B$ coming from the categories with an action of a group $\Gamma$ considered in \cite{De}, there is a non-degenerate GSP realizing it. More precisely, the endomorphism ring of a $\Gamma$-stable cluster-tilting object in the stable category of a category constructed in \cite{De} can be realized by a non-degenerate GSP (it is the case because $\Gamma$-$2$-cycles do not appear after mutations). In particular, the only potential for an acyclic group species is non-degenerate.
\end{rem}

Another proposition linking rigid and non-degenerate potentials can be adapted from \citeb{proposition 8.1 and corollary 8.2}{DeWeZe08}:

\begin{prop}
 Every rigid globally free GSP $(A, S)$ is $2$-acyclic and, in this case, $S$ is non-degenerate.
\end{prop}

As in \cite[\S 8]{DeWeZe08}, there exist group species without rigid potentials. The techniques of \cite[\S 8]{DeWeZe08} work also in the context of this article.

\section{Decorated representations and their mutations} \label{mutrep}

The aim of this section is to adapt the results of \cite[\S 10]{DeWeZe08}. We suppose here that for any $i \in I$, the characteristic of $K$ does not divide the cardinal of $\Gamma_i$.

Following \citeb{definition 10.1}{DeWeZe08},

\begin{df}
 A \emph{decorated representation} of a GSP $(A, S)$ is a pair $(X, V)$ where $X$ is a $\Pr(A, S)$-module and $V$ is a $E$-module.
\end{df}

In the following, we will look at pairs consisting of a GSP $(A, S)$ and a decorated representation of it. We will denote this type of objects by $(A, S, X, V)$ and call them \emph{group species with potential and decorated representation} (GSPDR).

Following \citeb{definition 10.2}{DeWeZe08},

\begin{df}
 A right-equivalence between two GSPDRs $(A, S, X, V)$ and $(A', S', X', V')$ is a triple $(\phi, \psi, \eta)$ such that:
 \begin{itemize}
  \item $\phi: E \lla A \rra \rightarrow E \lla A' \rra$ is a right-equivalence from $(A, S)$ to $(A', S')$ (see definition \ref{reqg}); 
  \item $\psi: X \rightarrow X'$ is a linear isomorphism such that the following diagram commutes:
   $$\xymatrix{
    X \ar[r]^{u_X} \ar[d]_{\psi} & X \ar[d]^{\psi} \\
    X' \ar[r]_{\phi(u)_{X'}} & X'
   }$$
   for any $u \in E \lla A \rra$;
  \item $\eta: V \rightarrow V'$ is an isomorphism.
 \end{itemize}
\end{df}

Using proposition \ref{eqenltriv}, for each GSPDR $(A, S, X, V)$, the decorated representation $(X, V)$ can be seen as a representation of $(A_{\red}, S_{\red})$. Thus, we can call $(A_{\red}, S_{\red}, X, V)$ the \emph{reduced part} of $(A, S, X, V)$. As in \cite[proposition 10.5]{DeWeZe08}, the right-equivalence class of the reduced part of a GSPDR is fully determined by the right-equivalence class of this GSPDR.

Now, we can define the mutation of a GSPDR $(A, S, X, V)$. Let $k \in I$. Our aim is to define a GSPRD $\mu_k(A,S,X,V) = (A', S', X', V')$ such that $(A', S') = \mu_k(A, S)$. Denote:
 $$X_{\inW} = X \tens_E A E_k \quad \text{and} \quad X_{\out} = X \tens_E A^* E_k.$$
Thus, we can define two right $E_k$-module morphisms. One, $\alpha$, from $X_{\inW}$ to $X_k = X E_k$ which is the application $(x \tens a) \mapsto xa$ and one from $X_k$ to $X_{\out}$ which is defined by
 $$\beta(x) = \sum_{b \in \mathcal{B}} x b \tens b^*$$
 which does not depend on the basis $\mathcal{B}$ of $E_k A$. Observe also that we have a canonical sequence of isomorphisms:
 \begin{align*}
  \Hom_{E_k}(X_{\out}, X_{\inW}) & \simeq \Hom_E(X \tens_E A^* E_k \tens_{E_k} E_k A^*, X) \\
   & \simeq \Hom_E (X \tens_E (A E_k A)^*, X)
 \end{align*}
 It is not hard to see that $[x \tens \xi \mapsto x (\partial_\xi S)] \in \Hom_E (X \tens_E (A E_k A)^*, X)$. Let $\gamma$ be the corresponding element of $\Hom_{E_k}(X_{\out}, X_{\inW})$.

So we get, as in \cite{DeWeZe08} a commutative diagram of right $E_k$-modules:
 $$\xymatrix{
  & X_k \ar[dr]^\beta & \\
  X_{\inW} \ar[ur]^\alpha & & X_{\out} \ar[ll]^\gamma
 }$$
with $\alpha \gamma = \gamma \beta = 0$ \cite[lemma 10.6]{DeWeZe08}. For $i \in I$, define:
 $$X'_i = \left\{\begin{array}{ll} 
  \displaystyle X_i & \quad \text{if } i \neq k \\
  \displaystyle \frac{\ker \gamma}{\ima \beta} \oplus \ima \gamma \oplus \frac{\ker \alpha}{\ima \gamma} \oplus V_i & \quad \text{if } i = k
       \end{array} \right. $$
  and
  $$V'_i = \left\{\begin{array}{ll} 
   \displaystyle V_i & \quad \text{if } i \neq k \\
   \displaystyle \frac{\ker \beta}{\ker \beta \inter \ima \alpha} & \quad \text{if } i = k
       \end{array} \right.
 $$
 
To get the structure of an $\Pr(A', S')$-module on $X'$, we must define the way $\tilde A$ acts on it where $(\tilde A, \tilde S) = \tilde \mu_k(A, S)$ (as $\Pr(A', S') \simeq \Pr(\tilde A, \tilde S)$). Recall from \S \ref{mut}, that 
 $$\tilde A = \bar E_k A \bar E_k \oplus A E_k A \oplus (E_k A)^* \oplus (A E_k)^*.$$
First of all, $\bar E_k A \bar E_k \oplus A E_k A \subset \bar E_k E \lla A \rra \bar E_k$ and for the vertices outside $k$, $X'_k = X_k$. Therefore, we can take the same action for this part of $\tilde A$. For the rest, we have $\tilde A E_k = A^* E_k$ and $\tilde A^* E_k = A E_k$ and therefore, we have to define:
 $$\alpha': X'_{\inW} = X' \tens_E \tilde A E_k = X \tens_E A^* E_k = X_{\out} \rightarrow X'_k$$
and 
 $$\beta': X'_k \rightarrow X'_{\out} = X' \tens_E \tilde A^* E_k = X \tens_E A E_k = X_{\inW}$$

 As in \cite{DeWeZe08}, we have to choose a \emph{splitting data}:
 \begin{itemize}
  \item let $\rho: X_{\out} \surj \ker \gamma$ be a splitting of $\ker \gamma \inj X_{\out}$ in the category $\md E_k$ (it is possible, as the characteristic of $K$ does not divide the cardinal of $\Gamma_k$);
  \item let $\sigma: \ker \alpha/\ima \gamma \inj \ker \alpha$ a splitting of $\ker \alpha \surj \ker \alpha / \ima \gamma$ in $\md E_k$.
 \end{itemize}
 Now, using the direct sum decomposition 
 $$X'_k = \frac{\ker \gamma}{\ima \beta} \oplus \ima \gamma \oplus \frac{\ker \alpha}{\ima \gamma} \oplus V_i,$$
 define
 $$\alpha' = \left( \begin{matrix} - \pi \rho \\ - \gamma \\ 0 \\ 0 \end{matrix} \right) \quad \text{and} \quad \beta' = \left( \begin{matrix} 0 & \iota & \iota \sigma & 0 \end{matrix} \right)$$
 where $\pi$ designs the canonical projection and $\iota$ the canonical injections.

Again, \citeb{proposition 10.7}{DeWeZe08} can be adapted:

\begin{prop}
 The above definition gives rise to a decorated representation of $(\tilde A, \tilde S)$ and, therefore, through the isomorphism $\Pr(\tilde A, \tilde S) \simeq \Pr(A', S')$, to a decorated representation of $(A', S')$.
\end{prop}

\begin{nt}
 We denote 
 $$\tilde \mu_k(A, S, X, V) = (\tilde A, \tilde S, X', V') \quad \text{and} \quad \mu_k(A, S, X, V) = (A', S', X', V').$$
\end{nt}

We can adapt \citeb{proposition 10.9}{DeWeZe08}:

\begin{prop}
 The isomorphism class of the GSPDR $\tilde \mu_k(A, S, X, V)$ does not depend on the choice of the splitting data.
\end{prop}

and \citeb{proposition 10.10 and corollary 10.12}{DeWeZe08}:

\begin{prop}
 The right-equivalence classes of the GSPDRs $$\tilde \mu_k(A, S, X, V) \quad \text{and} \quad \mu_k(A, S, X, V)$$ depend only on the right-equivalence class of $(A, S, X, V)$.
\end{prop}

Now an important theorem whose proof is the same as the one of \citeb{theorem 10.13}{DeWeZe08}:

\begin{thm}
 On the right-equivalence classes of GSPDRs which are $2$-acyclic at $k$, $\mu_k$ is an involution.
\end{thm}

It is easy to define the notion of a direct sum of two decorated representations of a GSP and, therefore, the notion of an indecomposable decorated representation of a GSP. Thus, as $\mu_k$ clearly commutes with this type of direct sums, $\mu_k$ acts on GSPs with indecomposable decorated representations. We call a GSPDR $(A, S, X, V)$ positive if $V = \{0\}$ and negative if $X = \{0\}$. Moreover, it is called \emph{simple} at $i \in I$ if $X \oplus V$ is an indecomposable $E_i$-module. Then we adapt \citeb{proposition 10.15}{DeWeZe08}:

\begin{prop} \label{indecrd}
 An indecomposable GSPDR is either positive, or negative simple. The mutation $\mu_k$ exchange a positive simple at $k$ with the corresponding negative simple at $k$. Moreover, it is the only case where a mutation interchanges positive and negative indecomposable GSPDRs.
\end{prop}

As in \cite[\S 6]{DeWeZe}, denote, for $k \in I$ and $X, X' \in \md \Pr(A, S)$,
 $$\Hom^{[k]}_{\Pr(A, S)}(X, X') = \left\{f \in \Hom_{\Pr(A, S)}(X, X') \,|\, f\restr{X \bar E_k} = 0 \right\}.$$

Cite now easy to adapt \citeb{proposition 6.1}{DeWeZe}:

\begin{prop} \label{modconf}
 The mutation $\mu_k$ induces an isomorphism
 $$\frac{\Hom_{\Pr(A, S)}(X, X')}{\Hom^{[k]}_{\Pr(A, S)}(X, X')} \simeq \frac{\Hom_{\Pr(\mu_k(A, S))}(\mu_k(X), \mu_k(X'))}{\Hom^{[k]}_{\Pr(\mu_k(A, S))}(\mu_k(X), \mu_k(X'))}.$$
\end{prop}

\begin{rem}
 As claimed in \cite[\S 6]{DeWeZe}, the isomorphism of proposition \ref{modconf} can be seen as a functorial isomorphism by introducing adapted quotient categories.
\end{rem}

\section{$F$-polynomials and $\g$-vectors of decorated representations} \label{fpol}

The aim of this section is to define the notions of the $F$-polynomial and the $\g$-vector of a GSPDR and to give a link with the usual notion (see \cite{FoZe07}). It is an extension of \cite{DeWeZe}. As before, $(I, (\Gamma_i))$ and therefore $E$ are fixed. We suppose also that the characteristic of $K$ does not divide any of the cardinals of the groups $\Gamma_i$. We suppose moreover that $K$ is algebraically closed and that all the $\Gamma_i$ are commutative (as seen in section \ref{exc}, this case is sufficient to realize skew-symmetrizable exchange matrices).

\begin{nt}
 For any $i \in I$, denote $\irr_i = \irr(\Gamma_i)$ the set of isomorphism classes of irreducible representations of $\Gamma_i$. One defines $\irr = \bigcup_{i \in I} \{i\} \times \irr_i$ and for $i \in I$, $C_i = K_0(\Gamma_i) \simeq \Z^{\irr_i}$. We also denote $C = K_0(E) = \bigoplus_{i \in I} C_i \simeq \Z^{\irr}$. If $V \in \md E$ (resp. $V \in \md E_i$), $[V]$ is its class in $C$ (resp. in $C_i$). If $\egr \in C$ (resp. $\egr \in C_i$) and $(j, \rho) \in \irr$ (resp. $\rho \in \irr_i$) then $\egr_{j,\rho}$ (resp. $\egr_\rho$) is the coefficient of $(j, \rho)$ (resp. $\rho$) in $\egr$.

 If $(Y_j)_{j \in \irr}$ (resp. $(Y_j)_{j \in \irr_i}$) is a family of indeterminates or of elements of a ring, and $\egr \in C$ (resp. $\egr \in C_i$), one denotes
 $$Y^{\egr} = \mathop{\prod_{j \in \irr}}_{(\text{resp. } j \in \irr_i)} Y_j^{\egr_j}.$$

 If $(A, S)$ is a GSP, $X$ a representation of it, $[X]$ is its class, seen as an $E$-module, in $C$. If $\egr \in C$ then $\Gr_{\egr}(X)$ is the Grassmanian of the $\Pr(A, S)$-submodules $X'$ of $X$ such that $[X'] = \egr$.
\end{nt}

Let $(A, S, X, V)$ be a GSPDR, we recall the diagram of section \ref{mutrep}, by changing a little the notation:
$$\xymatrix{
  & X(k) \ar[dr]^{\beta_k} & \\
  X_{\inW}(k) \ar[ur]^{\alpha_k} & & X_{\out}(k) \ar[ll]^{\gamma_k}
 }$$

\begin{df}
 One defines the $F$-polynomial $F_X$ of $X$ to be a polynomial in $\Z\left[(Y_i)_{i \in \irr}\right]$ defined by:
 $$F_X(Y) = \sum_{\egr \in C} \chi\left(\Gr_{\egr}(X)\right) Y^{\egr}$$
 where $\chi$ is the Euler characteristic. One define also the $\g$-vector $\g_{X,V} = \left(g_k\right)_{k \in I} \in C = \bigoplus_{k \in I} C_k$ by
 $$g_k = [\ker \gamma_k] - [X(k)] + [V(k)].$$
 With the same indexing, define $\h_{X,V} = \left(h_k\right)_{k \in I}$ by
 $$h_k = -[\ker \beta_k].$$ 
\end{df}

\begin{nt}
 If $(Y)$ is a family of indeterminates, we denote by $\Q_+(Y)$ the free commutative semifield generated by its elements. If $(y)$ is a family of elements of a commutative semifield, we denote by $\Q_+(y)$ the subsemifield generated by its elements.
\end{nt}

Then, it is easy to adapt \citeb{proposition 3.1}{DeWeZe}, \citeb{proposition 3.2}{DeWeZe} and \citeb{proposition 3.3}{DeWeZe}:

\begin{prop}
 The polynomial $F_X(Y)$ has constant term $1$ and maximum term (for divisibility of monomials) $Y^{[X]}$.
\end{prop}

\begin{prop}
 If $X'$ is another $\Pr(A, S)$-module then $F_{X \oplus X'} = F_X F_{X'}$. 
\end{prop}

\begin{prop} \label{lienFh}
 If $F_X \in \Q_+(Y)$, then $F_X$ can by evaluated in the semifield $\Trop(Y')$ where $(Y')_{i \in \irr}$ is a family of indeterminates. Then $\h_X$ and $F_X$ are related by the following formula:
 $$Y'^{\h_X} = F_X\restr{\Trop(Y')} \left(Y_{i,\rho}'^{-1} Y'^{[\rho \tens_{E_i} E_i A^*]}\right)_{(i, \rho) \in \irr}.$$ 
\end{prop}

\begin{demo}
 We follow the proof of \cite{DeWeZe}. Remark that for any $\egr \in C$,
 $$\left(Y^{\egr}\right)\restr{\Trop(Y')} \left(Y_{i,\rho}'^{-1} Y'^{[\rho \tens_{E_i} E_i A^*]}\right)_{(i, \rho) \in \irr} = Y'^{-\egr + [\egr \tens_E A^*]}.$$
 For $i \in I$, the exponent of $Y'_i = \left(Y_{i,\rho}\right)_{\rho \in \irr_i}$ can be rewritten as
 $$-\egr_i + [\egr \tens_E A^* E_i]$$
 which can be interpreted as
 $$-[X'(i)] + [X'_{\out}(i)]$$
 for any submodule $X'$ of $X$ such that $[X'] = \egr$. Thus, the end of the proof is the same as in \cite{DeWeZe}. \cqfd
\end{demo}

Recall the definition of a $Y$-seed:

\begin{df}[\citeb{\S 2}{DeWeZe}]
 A \emph{$Y$-seed} is a pair $(y, B)$ where $y$ is a family of elements of a semifield indexed by $I$ and $B$ is a skew-symmetrizable matrix. For $k \in I$, we define $\mu_k(y, B) = (y', \mu_k(B))$ where, for $i \in I$, 
 $$y'_i = \left\{\begin{array}{ll}
   y_i^{-1} & \quad \text{if } i = k \\
   y_i y_k^{\max(0, b_{ki})} (1+y_k)^{-b_{ki}} & \quad \text{if } i \neq k. \end{array}\right.$$   
\end{df}

Now, define the notion of an extended $Y$-seed:

\begin{df}
 A \emph{extended $Y$-seed} is a pair $(y, (A, S))$ where $y$ is a family of elements of a semifield indexed by $\irr$ and $(A, S)$ is a non-degenerate GSP. For $k \in I$, we define $\mu_k(y, (A, S)) = (y', \mu_k(A, S))$ where, for $(i, \rho) \in \irr$, 
 $$y'_{i, \rho} = \left\{\begin{array}{ll}
   y_{i, \rho}^{-1} & \quad \text{if } i = k \\
   y_{i,\rho} y_k^{[\rho \tens_{E_i} A_{ik}]} (1+y_k)^{[\rho \tens_{E_i} A_{ki}^*] - [\rho \tens_{E_i} A_{ik}]} & \quad \text{if } i \neq k. \end{array}\right.$$
\end{df}

\begin{rem}
  The mutation of extended $Y$-seeds is involutive.
\end{rem}

\begin{df}
 A $Y$-seed or an extended $Y$-seed will be called \emph{free} if its variables $y$ are algebraically independent.
\end{df}

\begin{rem}
 The notion of freeness for a $Y$-seed (or an extended $Y$-seed) is stable under mutations. The semifield $\Z_+(y)$ and the algebra $\Z[y]$ are also stable under mutation, as the mutation is involutive.
\end{rem}

\begin{df}
 Let $(y, (A, S))$ be a free extended $Y$-seed and $(z, B(A))$ be a $Y$-seed (for the same $A$). The following morphism of algebra is called the \emph{specialization map}:
 \begin{align*}
  \Phi_{y \rightarrow z}: \Z_+(y) & \rightarrow \Z_+(z) \\
  y_{i, \rho} & \mapsto z_i.
 \end{align*}
 The analogous for $\Z[y]$ and $\Z[z]$ is also denoted by $\Phi$.
\end{df}

\begin{prop} \label{mutPhi}
 Let $(y, (A, S))$ be a free extended $Y$-seed such that $(A, S)$ is a locally free GSP, and $(z, B(A))$ be a $Y$-seed. Let $k \in I$. Denote $y' = \mu_k(y)$, and $z' = \mu_k(z)$. Then, $\Phi_{y' \rightarrow z'} = \Phi_{y \rightarrow z}$.
\end{prop}

\begin{demo}
 As $y'$ generates $\Z_+(y') = \Z_+(y)$, it is enough to look at this for the $y'_{i, \rho}$ for $(i, \rho) \in \irr$. By definition,
 $$\Phi_{y'->z'}\left(y'_{i, \rho}\right) = {z'_i}$$
 If $i = k$, then
 $$\Phi_{y->z}\left(y'_{i, \rho}\right) = \Phi_{y->z}\left(y_{i, \rho}^{-1}\right) = z_i^{-1} = {z'_i}.$$
 If $i \neq k$, then
 \begin{align*}
  \Phi_{y->z}\left(y'_{i, \rho}\right) &= \Phi_{y->z}\left( y_{i,\rho} y_k^{[\rho \tens_{E_i} A_{ik}]} (1+y_k)^{[\rho \tens_{E_i} A_{ki}^*] - [\rho \tens_{E_i} A_{ik}]} \right) \\
   &= z_i \prod_{\sigma \in C_k} \left[ z_k^{[\rho \tens_{E_i} A_{ik}]_\sigma} (1+z_k)^{[\rho \tens_{E_i} A_{ki}^*]_\sigma - [\rho \tens_{E_i} A_{ik}]_\sigma} \right] \\
   &= z_i \left[ z_k^{\dim_K (\rho \tens_{E_i} A_{ik})} (1+z_k)^{\dim_K (\rho \tens_{E_i} A_{ki}^*) - \dim_K(\rho \tens_{E_i} A_{ik})} \right] \\
   &= z_i \left[ z_k^{\dim_{E_i} A_{ik}} (1+z_k)^{\dim_{E_i} A_{ki}^* - \dim_{E_i} A_{ik}} \right] \\
   &= z_i \left[ z_k^{\max(0, b_{ki})} (1+z_k)^{-b_{ki}} \right] = z_i'
 \end{align*}
 (here we use the fact that every considered irreducible representation is of dimension $1$, as the considered groups are commutative and $K$ is algebraically closed). \cqfd
\end{demo}

To make the relation with $F$-polynomials and $\g$-vectors in cluster algebras, we need the following adaptation of \cite[lemma 5.2]{DeWeZe}:

\begin{prop} \label{propF}
 Let $(A, S, X, V)$ be a GSPDR such that $(A, S)$ is non-degenerate. Let $k \in I$. Denote $(A', S', X', V') = \mu_k(A, S, X, V)$. Suppose also that the extended $Y$-seed $(y', (A',S'))$ is obtained from $(y, (A, S))$ by the mutation at $k$. Denote $\g_{X, V} = (g_i)_{i \in I}$, $\g_{X', V'} = (g'_i)_{i \in I}$, $\h_{X, V} = (h_i)_{i \in I}$ and $\h_{X', V'} = (h'_i)_{i \in I}$. Then
 \begin{enumerate} 
  \item $\g_{X, V} = \h_{X, V} - \h_{X', V'}$;
  \item one has
   $$(y_k + 1)^{h_k} F_X(y) = (y'_k + 1)^{h'_k} F_{X'}(y')$$
   where 
   $$(y_k + 1)^{h_k} = \prod_{i \in \irr_k} (y_{(k,i)} + 1)^{h_{ki}};$$
  \item for any $j \in I$,
    $$g'_j = \left\{ \begin{array}{ll}
      -g_j & \quad \text{if } j = k \\
       g_j + \left[g_k \tens_{E_k} A_{kj}\right] - \left[h_k \tens_{E_k} A_{kj}\right] + \left[h_k \tens_{E_k} A_{jk}^*\right]  & \quad \text{if } j \neq k.
     \end{array} \right.$$
 \end{enumerate}
\end{prop}

\begin{demo}
 \begin{enumerate}
  \item By definition, for $i \in I$, $g_i = [\ker \gamma_i] - [X(i)] + [V(i)]$, $h_i = -[\ker \beta_i]$ and $h'_i = -[\ker \beta'_i]$ (where $\beta'$ is the analogous of $\beta$ for $(X', V')$). So it is enough to prove that
  $$[\ker \gamma_i] + [V_i] + [\ker \beta_i] = [X(i)] + [\ker \beta'_i].$$
  From the definition of $\beta'_i$ given in section \ref{mutrep}, it is easy to see that $\ker \beta'_i \simeq \ker (\gamma_i) / \ima (\beta_i) \oplus V_i$. And, therefore, the searched equality reduces to 
  $$[\ima \beta_i] + [\ker \beta_i] = [X(i)]$$
  which is obvious.
  \item We follow the proof of \cite[lemma 5.2]{DeWeZe}. Let $\egr \in C$ and $\egr'$ its projection in $\bigoplus_{i \neq k} C_i$. Let $X_0 = X \bar E_k$ which is a $\bar E_k \Pr(A, S) \bar E_k$-module. For any $\bar E_k \Pr(A, S) \bar E_k$-submodule $W$ of $X_0$, one can define
  $$W_{\inW}(k) = W \tens_{\bar E_k} A E_k \subset X_{\inW}(k) \quad \text{and} \quad W_{\out}(k) = W \tens_{\bar E_k} A^* E_k \subset X_{\out}(k)$$
  which are well defined because $(A, S)$ has no loop (and therefore $X_{\inW} = X \tens_{\bar E_k} A E_k$ and $X_{\out} = X \tens_{\bar E_k} A^* E_k$).

 For $\rgr, \sgr \in C_k$, define $Z_{\egr', \rgr, \sgr}(X)$ to be the subvariety of $\Gr_{\egr'}(X_0)$ consisting of the $W$ satisfying
 \begin{itemize}
  \item $\left[\alpha_k\left(W_{\inW}(k)\right)\right] = \rgr$;
  \item $\left[\beta_k^{-1}\left(W_{\out}(k)\right)\right] = \sgr$;
  \item $\alpha_k\left(W_{\inW}(k)\right) \subset \beta_k^{-1}\left(W_{\out}(k)\right)$.
 \end{itemize}
 Define also the variety
 $$\tilde Z_{\egr, \rgr, \sgr}(X) = \left\{W \in \Gr_{\egr}(X) \,|\, W \bar E_k \in Z_{\egr', \rgr, \sgr}(X) \right\}$$
 so that, by an easy computation, $\tilde Z_{\egr, \rgr, \sgr}(X)$ is a fiber bundle over $Z_{\egr', \rgr, \sgr}(X)$ with fiber $\Gr_{e_k - \rgr}(\sgr - \rgr)$ (where, by abuse of notation, we identify $\sgr - \rgr \geq 0 $ with any of its representatives in $\md E_k$, and $\Gr_{e_k - \rgr}(\sgr - \rgr) = \emptyset$ if $e_k - \rgr$ or $\sgr-\rgr$ are not nonnegative). Hence, using the easy fact that $\Gr_{\egr}(X)$ is the disjoint union of the $\tilde Z_{\egr, \rgr, \sgr}(X)$, we obtain, as every considered irreducible representation is of dimension $1$,
 $$\chi\left(\Gr_{\egr}(X)\right) = \sum_{\rgr, \sgr \in C_k} \vecv{\sgr - \rgr}{e_k - \rgr} \chi\left(Z_{\egr', \rgr, \sgr}(X)\right).$$
 where, for any $\rgr_1, \rgr_2 \in C_k$, 
 $$\vecv{\rgr_1}{\rgr_2} = \prod_{\rho \in \ind_k} \vecv{\rgr_{1, \rho}}{\rgr_{2,\rho}}$$
 Then, substituting this expression in the definition of $F_X$, we obtain:
 \begin{align*}
  F_X(y) &= \sum_{\egr \in C} \left[ \sum_{\rgr, \sgr \in C_k} \vecv{\sgr - \rgr}{e_k - \rgr} \chi\left(Z_{\egr', \rgr, \sgr}(X)\right) \right] y^{\egr} \\
   &= \mathop{\sum_{\egr' \in \bigoplus_{i \neq k} C_i}}_{\rgr, \sgr \in C_k} \chi\left(Z_{\egr', \rgr, \sgr}(X)\right) y^{\egr'}  \sum_{e_k \in C_k} \vecv{\sgr - \rgr}{e_k - \rgr} y_k^{e_k} \\
   &= \mathop{\sum_{\egr' \in \bigoplus_{i \neq k} C_i}}_{\rgr, \sgr \in C_k} \chi\left(Z_{\egr', \rgr, \sgr}(X)\right) y^{\egr' + \rgr}  (1 + y_k)^{\sgr - \rgr}.
 \end{align*}
 
 Now, as in \cite{DeWeZe}, we have easily that
 $$Z_{\egr', \rgr, \sgr}(X) = Z_{\egr', \bar \rgr, \bar \sgr}(X')$$
 where
 $$\bar \rgr = \left[\egr' \tens_{\bar E_k} A^* E_k \right] - h_k - \sgr \quad \text{and} \quad \bar \sgr = \left[\egr' \tens_{\bar E_k} A E_k \right] - h'_k - \rgr.$$
 Using this, one gets
 \begin{align*}
  (1+y'_k)^{h'_k} F_{X'}(y') &= \mathop{\sum_{\egr' \in \bigoplus_{i \neq k} C_i}}_{\bar \rgr, \bar \sgr \in C_k} \chi\left(Z_{\egr', \bar \rgr, \bar \sgr}(X')\right) y'^{\egr' + \bar \rgr}  (1 + y'_k)^{h'_k + \bar \sgr - \bar \rgr} \\
   &= \mathop{\sum_{\egr' \in \bigoplus_{i \neq k} C_i}}_{\rgr, \sgr \in C_k} \chi\left(Z_{\egr', \rgr, \sgr}(X)\right) y'^{\egr'} y_k^{- \bar \sgr - h'_k}  (1 + y_k)^{h'_k + \bar \sgr - \bar \rgr} \\
   &= \mathop{\sum_{\egr' \in \bigoplus_{i \neq k} C_i}}_{\rgr, \sgr \in C_k} \chi\left(Z_{\egr', \rgr, \sgr}(X)\right) y^{\egr' + \rgr}  (1 + y_k)^{h_k + \sgr - \rgr} \\
   &= (1+y_k)^{h_k} F_{X}(y)
 \end{align*}
 \item As $g_k = h_k - h'_k$, $g'_k = - g_k$. If $j \neq k$, the equality we want to prove becomes, using again $g_k = h_k - h'_k$,
   $$[\ker \gamma'_j] - \left[\ker \beta'_k \tens_{E_k} A_{kj}\right]  = [\ker \gamma_j] - \left[\ker \beta_k \tens_{E_k} A_{jk}^*\right]$$
  and, up to a possible exchange of $(A, S, X, V)$ and $(A', S', X', V')$, we can suppose that $A_{kj} = 0$ (because $A$ is $2$-acyclic) and therefore, we have to prove that
  $$[\ker \gamma'_j]  = [\ker \gamma_j] - \left[\ker \beta_k \tens_{E_k} A_{jk}^*\right].$$ 
 Let 
  $$(\tilde A, \tilde S, \tilde X, \tilde V) = \tilde \mu_k(A, S, X, V)$$
 in such a way that $(A', S')$ is right-equivalent to $(\tilde A, \tilde S)_{\rd}$. In this setting, one will prove that
  $$[\ker \tilde \gamma_j]  = [\ker \gamma_j] - \left[\ker \beta_k \tens_{E_k} A_{jk}^*\right].$$ 
 We can decompose
  $$X_{\out}(j) = X \tens_E A^* E_j = X(k) \tens_{E_k} A^*_{jk} \oplus X \bar E_k \tens_{\bar E_k} \bar E_k A^* E_j$$
 and we get
  $$\tilde X_{\out}(j) = X_{\out}(k) \tens_{E_k} A^*_{jk} \oplus X \bar E_k \tens_{\bar E_k} \bar E_k A^* E_j$$
 and
  $$\tilde X_{\inW}(j) = \tilde X(k) \tens_{E_k} \tilde A_{kj} \oplus X_{\inW}(j) = X'(k) \tens_{E_k} A^*_{jk} \oplus X_{\inW}(j).$$
 Along these decompositions, one has:
  $$\gamma_j = \vech{\psi \circ (\beta_k \tens_{E_k} A^*_{jk})}{\eta} \quad \text{and} \quad \tilde \gamma_j = \mat{\alpha'_k \tens_{E_k} A^*_{jk}}{0}{\psi}{\eta}$$
 where $\psi: X_{\out}(k) \tens_{E_k} A^*_{jk} \rightarrow X_{\inW}(j)$ and $\eta: X \bar E_k \tens_{\bar E_k} \bar E_k A^* E_j \rightarrow X_{\inW}(j)$ are two $E_j$-modules morphisms (basically speaking, these two morphisms encode the part of $\gamma_j$ which is not modified by the mutation at $k$). Using definitions of section \ref{mutrep}, we get easily that $\ker \alpha'_k = \ima \beta_k$ and we get an exact sequence of $E_j$-modules:
 $$0 \rightarrow \ker \beta_k \tens_{E_k} A_{jk}^* \oplus \{0\} \rightarrow \ker \gamma_i \xrightarrow{f} \ker \tilde \gamma_i \rightarrow 0$$
 where, along the previous decompositions $$f(u, v) = ((\beta_k \tens_{E_k} A^*_{jk}) u, v).$$
 This short exact sequence implies that
  $$[\ker \tilde \gamma_j]  = [\ker \gamma_j] - \left[\ker \beta_k \tens_{E_k} A_{jk}^*\right].$$ 
  To finish, it remains to prove that $[\ker \tilde \gamma_j] = [\ker \gamma'_j]$. The proof is the same than in \cite{DeWeZe}. \cqfd
 \end{enumerate}
\end{demo}

\begin{df}
 For any GSPDR $(A, S, X, V)$, we define in the following way the reduced $\g$-vectors, $\h$-vectors and $F$-polynomials:
 \begin{itemize}
  \item for $i \in I$, let $\check \g_{X, V} = (\check g_i)_{i \in I}$ defined by $\check g_i = \dim_K g_i$ where $(g_i)_{i \in I} = \g_{X, V}$;
  \item for $i \in I$, let $\check \h_{X, V} = (\check h_i)_{i \in I}$ defined by $\check h_i = \dim_K h_i$ where $(h_i)_{i \in I} = \h_{X, V}$;
  \item $\check F_X = \Phi_{Y \rightarrow Z} (F_X)$ where $(Y_i)_{i \in \irr}$ and $(Z_i)_{i \in I}$ are families of indeterminates.
 \end{itemize}
\end{df}

\begin{cor} \label{vermat}
 Let $(A, S, X, V)$ be a GSPDR such that $(A, S)$ is non-degenerate and locally free. Let $k \in I$. Denote $$(A', S', X', V') = \mu_k(A, S, X, V).$$ Suppose also that the $Y$-seed $(z', B(A'))$ is obtained from $(z, B(A))$ by the mutation at $k$. Denote $\check \g_{X, V} = (\check g_i)_{i \in I}$, $\check \g_{X', V'} = (\check g'_i)_{i \in I}$, $\check \h_{X, V} = (\check h_i)_{i \in I}$ and $\check \h_{X', V'} = (\check h'_i)_{i \in I}$. We also denote by $(b_{ij})_{i,j \in I}$ the coefficients of $B(A)$. Then
 \begin{enumerate} 
  \item \label{coro1} $\forall i \in I, \check g_i = \check h_i - \check h'_i$;
  \item \label{coro2} one has
   $$(z_k + 1)^{\check h_k} \check F_X(z) = (z'_k + 1)^{\check h'_k} \check F_{X'}(z');$$
  \item \label{coro3} for any $j \in I$,
    $$\check g'_j = \left\{ \begin{array}{ll}
      -\check g_j & \quad \text{if } j = k \\
       \check g_j + \max(0, b_{jk}) \check g_k - b_{jk} \check h_k  & \quad \text{if } j \neq k;
     \end{array} \right.$$
  \item \label{coro4} if $F_X \in \Q_+(Y)$, then $\check F_X \in \Q_+(Z)$. Then $\check \h_X$ and $\check F_X$ are related by the following formula:
 $$Z_0^{\check \h_X} = \check F_X\restr{\Trop(Z_0)} \left(Z_{0,i}^{-1} \prod_{j \neq i} Z_{0,j}^{\max(0,-b_{ji})}\right)_{i \in I}.$$ 
 \end{enumerate}
\end{cor}

\begin{demo}
 The points (\ref{coro1}) and (\ref{coro3}) are immediate consequences of proposition \ref{propF}. To prove (\ref{coro2}), it is enough to apply $\Phi_{y \rightarrow z}$ to the analogous identity in proposition \ref{propF} (for any extended free $Y$-seed $(y, (A, S))$) and then apply proposition \ref{mutPhi}. For (\ref{coro4}), remark that for any $(i, \rho) \in \irr$, 
 $$\Phi_{Y_0 \rightarrow Z_0}\left(Y_{0,i,\rho}^{-1} Y_0^{[\rho \tens_{E_i} E_i A^*]}\right) = Z_{0,i}^{-1} \prod_{j \neq i} Z_{0,j}^{\max(0,-b_{ji})}$$
 is independent of $\rho$ and therefore, it is easy to see that 
 \begin{align*} & \check F_X\restr{\Trop(Z_0)} \left(Z_{0,i}^{-1} \prod_{j \neq i} Z_{0,j}^{\max(0,-b_{ji})}\right)_{i \in I} \\ =& \Phi_{Y_0 \rightarrow Z_0}\left(F_X\restr{\Trop(Y_0)} \left(Y_{0,i,\rho}^{-1}  Y_0^{[\rho \tens_{E_i} E_i A^*]} \right)_{(i,\rho) \in I}\right) \\=& \Phi_{Y_0 \rightarrow Z_0} \left(Y_0^{\h_X}\right) = Z_0^{\check \h_X} \end{align*}
 using proposition \ref{lienFh}. \cqfd
\end{demo}

In \cite{FoZe07}, (see also \cite[\S 2]{DeWeZe}), Fomin and Zelevinsky defined the notions of the $F$-polynomials and the $\g$-vectors associated to a sequence of mutation. More precisely, for a skew-symmetrizable matrix $B$ (which will play the role of an initial seed), a sequence of indices $\ig = (i_1, i_2, \dots, i_n) \in I^n$ and $k \in I$, they define a polynomial $F^B_{k;\ig} \in \Z[Z_i]_{i \in I}$ and a vector $\g^B_{k;\ig} \in \Z^I$. 

\begin{df} \label{intM}
 Let $(A, S)$ be a non-degenerate GSP and $\ig = (i_1, \dots, i_n)$ be in $I^n$ and $V$ an $E$-module. We denote
 $$\left( A^{A,S}_{V; \ig}, S^{A,S}_{V; \ig}, X^{A,S}_{V; \ig}, V^{A,S}_{V; \ig} \right) = \mu_{i_1} \mu_{i_2} \dots \mu_{i_n} \left(\mu_{i_n} \dots \mu_{i_2} \mu_{i_1}(A,S), 0, V \right).$$
 Remark that $\left( A^{A,S}_{V; \ig}, S^{A,S}_{V; \ig}\right)$ is right-equivalent to $(A, S)$.
\end{df}

Thus, we can adapt theorem \cite[theorem 5.1]{DeWeZe}:

\begin{thm} \label{realFg}
 Let $(A, S)$ be a non-degenerate locally free GSP. Let $\ig = (i_1, i_2, \dots, i_n) \in I^n$, $k \in I$ and $\rho \in \irr_k$. Then
 $$\g^{B(A)}_{k;\ig} = \check \g_{X^{A, S}_{\rho; \ig}, V^{A, S}_{\rho; \ig}} \quad \text{and} \quad F^{B(A)}_{k;\ig} = \check F_{X^{A, S}_{\rho; \ig}}.$$
\end{thm}

\begin{demo}
 With corollary \ref{vermat}, it is the same proof as in \cite{DeWeZe}. \cqfd
\end{demo}

We get also this following, analogous to \citeb{corollary 5.3}{DeWeZe}:

\begin{cor}
 In the situation of theorem \ref{realFg}, suppose that $F^{B(A)}_{k; \ig} \neq 1$, hence $X^{A, S}_{\rho; \ig} \neq \{0\}$ and $V^{A, S}_{\rho; \ig} = \{0\}$ (see proposition \ref{indecrd}). Let $x^{B(A)}_{k; \ig}$ be the corresponding cluster variable in the coefficient-free cluster algebra. In other terms
 $$\left(\left(x^{B(A)}_{i; \ig}\right)_{i \in I}, B' \right) = \mu_{i_n} \dots \mu_{i_2} \mu_{i_1}\left(\left(x_i\right)_{i \in I} , B(A) \right).$$
 Then we have the following cluster character formula:
 $$x^{B(A)}_{k; \ig} = \prod_{i \in I} x_i^{-d_i} \sum_{\egr \in C} \chi\left(\Gr_{\egr}(X) \right) \prod_{i \in I} x_i^{- \rg \gamma_i + \sum_{j \in I} \left(\max(0, b_{ij}) e_j + \max(0, -b_{ij})(d_j - e_j) \right) }$$
 where $X = X^{A, S}_{\rho; \ig}$, $d_i = \dim_K X(i)$ and $e_i = \dim_K \egr_i$.
\end{cor}

\section{$\mathcal{E}$-invariant}

The aim of this part is analogous to \cite[\S 7, \S 8]{DeWeZe}. Let $(A, S, X, V)$ and $(A, S, X', V')$ be two GSPDRs with the same non-degenerate GSP. We denote:
 $$\langle X, X' \rangle = \dim_K \Hom_{\Pr(A, S)}(X, X').$$
 Define the three following integer functions:
 $$\Er^{\injW}(X, V; X', V') = \langle X, X' \rangle + \left([X] | \g_{X',V'}\right)$$
 $$\Er^{\sym}(X, V; X', V') = \Er^{\injW}(X,V;X',V') + \Er^{\injW}(X',V';X,V)$$
 $$\Er(X, V) = \Er^{\injW}(X,V;X,V)  = \frac{\Er^{\sym}(X,V;X,V)}{2}$$
 where $[X] \in C$ is the class of $X$ seen as an $E$-module, and, for $\egr, \egr' \in C$ (resp. $\egr, \egr' \in C_k$ for $k \in I$), 
 $$\left(\egr | \egr'\right) = \mathop{\sum_{i \in \irr}}_{(\text{resp. } i \in \irr_k)} \egr_i \egr'_i.$$
 
Then, we get, with the same proof as \citeb{theorem 7.1}{DeWeZe}:

\begin{thm}
 We have, for any $k \in I$, 
 \begin{align*} &\Er^{\injW}\left( \mu_k(X, V); \mu_k(X', V')\right) - \Er^{\injW}\left( X, V; X', V'\right) \\=& \left(\h_{\mu_k(X,V),k}|\h_{X',V',k}\right) - \left(\h_{X,V,k}|\h_{\mu_k(X',V'),k}\right).\end{align*}
 In particular, $\Er^{\sym}$ and $\Er$ are stable under mutations.
\end{thm}

\begin{demo}
 The only difference with \cite{DeWeZe} is that computations have to be done in the Grothendieck groups. Moreover, we have to worry about the skew-symmetrizability: with our convention, informally speaking, all $b_{ik}$ should be replaced by $-b_{ki}$ in the proof of \cite{DeWeZe}). For example, $$\sum_{i \in I} \max(0, b_{ik}) \dim_K X(i)$$ in \cite{DeWeZe} will be replaced here by $[X \tens_E A^* E_k]$ whose dimension is $$\sum_{i \in I} \max(0, -b_{ki}) \dim_K X(i)$$ if the GSP is locally free and $B = B(A)$. \cqfd
\end{demo}

We get also the following analogous of \citeb{corollary 7.2}{DeWeZe}:

\begin{cor}
 If $(X, V)$ is obtained by a sequence of mutations from a negative decorated representation $(\{0\}, V)$ then $\Er(X, V) = 0$.
\end{cor}

We denote by $A^{\op}$ the $(E, E)$-bimodule whose underlying vector space is $A$ and whose bimodule structure is given by $g \cdot a^{\op} \cdot h = (h^{-1} \cdot a \cdot g^{-1})^{\op}$ if $g \in \Gamma_i$ and $h \in \Gamma_j$ for some $i,j \in I$ and $\op : A \rightarrow A^{\op}$ comes from the identity of $A$. It is then easy to extend $\op$ to an anti-isomorphism of algebras $E \lla A \rra \rightarrow E \lla A^{\op} \rra$. Thus, $(X^*, V^*)$ is a decorated representation of the GSP $(A^{\op}, S^{\op})$ on the ring $E$, where for each $i \in I$, $X_i^*$ is contragredient to $X_i$, $V_i^*$ is contragredient to $V_i$ and $a^{\op}$ acts on $X^*$ as the transpose of $a$ for every $a \in A$. Thus, one gets the analogous of \citeb{proposition 7.3}{DeWeZe}:

\begin{prop}
 We have $\Er(X^*, V^*) = \Er(X, V)$.
\end{prop}

\begin{demo}
 As for any $i \in I$, the characteristic of $K$ does not divide $\# \Gamma_i$, we have an isomorphism of right $E$-modules
 \begin{align*}
  \left( X \tens_E A\right)^* &\rightarrow X^* \tens_E A^{*\op} \simeq X^* \tens_E A^{\op*} \\
  f & \mapsto \mathop{\sum_{x \in \mathcal{B}_X}}_{a \in \mathcal{B}_A} f(x \tens a) x^* \tens a^{*\op} \\
  \left(x \tens a \mapsto \sum_{i \in I} \sum_{g \in \Gamma_i} \frac{\phi(x g) \psi(g^{-1} a)}{\# \Gamma_i} \right)& \mapsfrom \phi \tens \psi^{\op}
 \end{align*}
 which does not depend of the bases $\mathcal{B}_X$ and $\mathcal{B}_A$ of $X$ and $A$. Thus, we have, as in \cite{DeWeZe},
 \begin{align*}
  \Er(X, V) =& \langle X, X \rangle + \left([X] | [X \tens_E A^*]\right) + \left([X] \vphantom{\left[\bigoplus_{i \in I} \ima \gamma_i\right]} \, \right| \left. [V]-[X]-\left[\bigoplus_{i \in I} \ima \gamma_i\right]\right) \\
  =& \langle X, X \rangle + \left([X \tens_E A] | [X]\right) + \left([X] \vphantom{\left[\bigoplus_{i \in I} \ima \gamma_i\right]} \, \right| \left. [V]-[X]-\left[\bigoplus_{i \in I} \ima \gamma_i\right]\right) \\
  =& \langle X^*, X^* \rangle + \left([(X \tens_E A)^*] | [X^*]\right) \\&+ \left([X^*] \vphantom{\left[\bigoplus_{i \in I} \ima \gamma^*_i\right]} \, \right| \left. [V^*]-[X^*]-\left[\bigoplus_{i \in I} \ima \gamma^*_i\right]\right) \\
  =& \langle X^*, X^* \rangle + \left([X^* \tens_E A^{\op *}] | [X^*]\right) \\&+ \left([X^*] \vphantom{\left[\bigoplus_{i \in I} \ima \gamma^*_i\right]} \, \right| \left. [V^*]-[X^*]-\left[\bigoplus_{i \in I} \ima \gamma^*_i\right]\right) \\
  =& \Er(X^*, V^*)
 \end{align*}
 where we used that \begin{align*} \left([X] | [X \tens_E A^*]\right) &= \dim_K \Hom_E (X, X \tens_E A^*) \\&= \dim_K \Hom_E (X \tens_E A, X) = \left([X \tens_E A] | [X]\right). \cqfd\end{align*}
\end{demo}

Hence, the following theorem has the same proof as \citeb{theorem 8.1}{DeWeZe} (note that all \cite[\S 10]{DeWeZe} can be easily adapted in this case):

\begin{thm}
 The $\Er$-invariant satisfies
 $$\Er(X, V) \geq \left(\left[\bigoplus_{i \in I} \ker \beta_i\right]\, \right|\left.\left[\bigoplus_{i \in I} \frac{\ker \gamma_i}{\ima \beta_i}\right] \right) + \left([X]|[V]\right).$$
\end{thm}

Then, we obtain the analogous of \citeb{corollary 8.3}{DeWeZe}:

\begin{cor} \label{eics}
 If $\Er(X, V) = 0$ then for each $(k, \rho) \in \irr$,
 \begin{enumerate}
  \item either $[M_k]_\rho = 0$ or $[V_k]_\rho = 0$;
  \item either $[\ker \gamma_k]_\rho = 0$ or $[\ker \gamma_k]_\rho = [\ima \beta_k]_\rho$.
 \end{enumerate}
\end{cor}

\section{Applications to cluster algebras} \label{clust}

We conclude here that the following conjectures of \cite{FoZe07} are true for skew-symmetrizable integer matrix which can be obtained from a non-degenerate GSP with abelian groups. In particular, every matrix of the form $DS$ where $D$ is diagonal with integer coefficients and $S$ is skew-symmetric with integer coefficients can be obtained in view of section \ref{nondeg}. Every exchange matrix corresponding to the situation described in \cite{De} (in particular every acyclic ones) can also be raised. Let $B$ be such a skew-symmetrizable integer matrix. We suppose moreover that some $(A, S)$ is fixed satisfying the hypothesis of section \ref{fpol} such that $B(A) = B$.

\begin{prop}[\citeb{conjecture 5.4}{FoZe07}]
 For any $\ig \in I^n$ and $k \in I$, $F^B_{k;\ig}$ has constant term $1$.
\end{prop}

\begin{prop}[\citeb{conjecture 5.5}{FoZe07}]
 For any $\ig \in I^n$ and $k \in I$, $F^B_{k;\ig}$ has a maximum monomial for divisibility order with coefficient $1$.
\end{prop}

These first two are immediate, as in \cite[\S 9]{DeWeZe}.

\begin{prop}[\citeb{conjecture 7.12}{FoZe07}]
 For any $\ig \in I^n$, $k \in I$, we denote by $k\ig$ the concatenation of $(k)$ and $\ig$. Let $j \in I$ and $(g_i)_{i \in I} = \g^B_{j; \ig}$ and $(g'_i)_{i \in I} = \g^{\mu_k(B)}_{j; k\ig}$. Then we have, for any $i \in I$,
 $$g'_i = \left\{\begin{array}{ll}
  -g_i & \quad \text{if } i = k; \\
  g_i + \max(0,b_{ik}) g_k - b_{jk} \min(g_k, 0) & \quad \text{if } i \neq k.
  \end{array} \right.$$
\end{prop}

\begin{demo}
 We need here to add some trick to the proof of \cite[\S 9]{DeWeZe}. Indeed, we need to prove, as in \cite{DeWeZe}, that
 $$\min(0,g_k) = h_k.$$
 But what we obtain by using corollary \ref{eics} is
 $$\min(0,g_{k,\rho}) = h_{k,\rho}$$
 for any $\rho \in \irr_k$. Moreover, we have, as seen before,
 $$g_k = \sum_{\rho \in \irr_k} g_{k, \rho} \quad \text{and} \quad h_k = \sum_{\rho \in \irr_k} h_{k, \rho}$$
 and therefore, what we need is equivalent to the fact that the $g_{k, \rho}$ are of the same sign. We will prove this with an indirect method. Retaining the notation of definition \ref{intM}, we get
 $$X_{E_j;\ig}^{A, S} = \sum_{\rho \in \irr_j} X_{\rho;\ig}^{A, S}$$
 and therefore, by linearity of $\g$, 
 $$\g_{X_{E_j;\ig}^{A, S}} = \sum_{\rho \in \irr_j} \g_{X_{\rho;\ig}^{A, S}}.$$
 Hence, we get:
 $$(\# \Gamma_j) g_k = \dim_K \left[\g_{X_{E_j;\ig}^{A, S}}\right]_k.$$
 In the same way,
 $$(\# \Gamma_j) h_k = \dim_K \left[\h_{X_{E_j;\ig}^{A, S}}\right]_k.$$
 Moreover, by an immediate induction using proposition \ref{propF}, as $[E_j]$ is the class of a free $E_j$-module, $\left[\g_{X_{E_j;\ig}^{A, S}}\right]_k$ and $\left[\h_{X_{E_j;\ig}^{A, S}}\right]_k$ are also free and therefore, their coefficients in term of the irreducible representations of $E_k$ are of the same sign. Hence, we obtain, by adding these components
 $$\min(0, (\# \Gamma_j) g_k) = (\# \Gamma_j) h_k$$
 and the rest follows as in \cite{DeWeZe}. Note that it implies also that the $g_{k, \rho}$ are of the same sign. \cqfd
\end{demo}

The three following propositions have the same proof than in \cite[\S 9]{DeWeZe}:

\begin{prop}[\citeb{conjecture 6.13}{FoZe07}]
 For any $\ig \in I^n$, the vectors $\g^B_{i;\ig}$ for $i \in I$ are sign-coherent. In other terms, for $i,i', j \in I$, the $j$-th components of $\g^B_{i;\ig}$ and $\g^B_{i';\ig}$ have the same sign.
\end{prop}

\begin{prop}[\citeb{conjecture 7.10(2)}{FoZe07}]
 For any $\ig \in I^n$, the vectors $\g^B_{i;\ig}$ for $i \in I$ form a $\Z$-basis of $\Z^I$.
\end{prop}

\begin{prop}[\citeb{conjecture 7.10(1)}{FoZe07}]
 For any $\ig, \ig' \in I^n$, if we have
 $$\sum_{i \in I} a_i \g_{i; \ig}^B = \sum_{i \in I} a'_i \g_{i; \ig'}^B$$
 for some nonnegative integers $(a_i)_{i \in I}$ and $(a'_i)_{i \in I}$, then there is a permutation $\sigma \in \mathfrak{S}_I$ such that
 for every $i \in I$, 
 $$a_i = a'_{\sigma(i)} \quad \text{and} \quad a_i \neq 0 \Rightarrow \g_{i;\ig}^B = \g_{\sigma(i);\ig'}^B \quad \text{and} \quad a_i \neq 0 \Rightarrow F_{i;\ig}^B = F_{\sigma(i);\ig'}^B.$$
 In particular, $F_{i; \ig}^B$ is determined by $\g_{i; \ig}^B$.
\end{prop}

\section{An example and a counterexample} \label{exemple}

\input{exemple.tex} 

\section*{Acknowledgments}

The author would like to thank Bernhard Keller, Sefi Ladkani, Bernard Leclerc and Dong Yang for discussions in relation to this subject.

\bibliographystyle{alphanum}
\bibliography{biblio}

\begin{thebibliography}{DWZ2}

\bibitem[Dem]{De}
Laurent Demonet.
\newblock Categorification of skew-symmetrizable cluster algebras.
\newblock arXiv: 0909.1633.

\bibitem[DR]{DlRi76}
Vlastimil Dlab and Claus~Michael Ringel.
\newblock Indecomposable representations of graphs and algebras.
\newblock {\em Mem. Amer. Math. Soc.}, 6(173):v+57, 1976.

\bibitem[DWZ1]{DeWeZe}
Harm Derksen, Jerzy Weyman, and Andrei Zelevinsky.
\newblock Quivers with potentials and their representations {II}: Applications
  to cluster algebras.
\newblock arXiv: 0904.0676.

\bibitem[DWZ2]{DeWeZe08}
Harm Derksen, Jerzy Weyman, and Andrei Zelevinsky.
\newblock Quivers with potentials and their representations. {I}. {M}utations.
\newblock {\em Selecta Math. (N.S.)}, 14(1):59--119, 2008.

\bibitem[FZ1]{FoZe02}
Sergey Fomin and Andrei Zelevinsky.
\newblock Cluster algebras. {I}. {F}oundations.
\newblock {\em J. Amer. Math. Soc.}, 15(2):497--529 (electronic), 2002.

\bibitem[FZ2]{FoZe07}
Sergey Fomin and Andrei Zelevinsky.
\newblock Cluster algebras. {IV}. {C}oefficients.
\newblock {\em Compos. Math.}, 143(1):112--164, 2007.

\end{thebibliography}

\end{document}